\documentstyle[12pt]{article}

  \newcommand{\const}{\rm const}
  \newcommand{\Var}{\rm Var}
  \newcommand{\supp}{\rm supp}
  \newcommand{\Law}{\rm Law}
  \newcommand{\mod}{\rm mod}
  \newcommand{\Dom}{\rm Dom}

\begin{document}

   \begin{center}

  {\bf   Prokhorov-Skorokhod continuity of random fields.  }\\

\vspace{4mm}

{\bf   A natural approach. }\\

\vspace{4mm}

   {\bf Ostrovsky E., Sirota L.}\\

\vspace{4mm}

 Israel,  Bar-Ilan University, department of Mathematic and Statistics, 59200, \\

\vspace{4mm}

E-mails: eugostrovsky@list.ru, \\
 sirota3@bezeqint.net \\

\vspace{5mm}

  {\bf Abstract} \\

 \end{center}

 \  We derive  in this article sufficient conditions in the natural terms for belonging of almost all the trajectories of the certain separable
 continuous in probability random field to the multivariate Prokhorov-Skorokhod space. \\
 \ We consider also as a consequence the Central Limit Theorem in this spaces. \\

\vspace{4mm}


\vspace{5mm}

\section{Definitions. Notations. Statement of problem.}

 \vspace{3mm}

 \  Let $  X = [0,1]^d, \ d = 1,2, 3, \ldots;   \ x = \vec{x} = (x_1,x_2, x_3, \ldots, x_d) \in X, \   $ and let $  \xi = \xi(x), \ x \in X $ be separable
stochastic continuous numerical valued  random process (r.p.)  or equally random field (r.f.,) in the general case. The correspondent probability space
will be denoted by $  (\Omega, B, {\bf P}) $ with expectation $  {\bf E} $ and variance $  \Var. $  Denote also by $ \ D(X) =  D[0,1]^d \ $ the multivariate
famous Prokhorov-Skorokhod space; we recall its definition further. \par

\vspace{4mm}

 \ {\bf  Our goal in this article is deriving some simple sufficient conditions for belonging of almost all trajectories of this random field to the
Prokhorov-Skorokhod space:  }

$$
{\bf P} (\xi(\cdot) \in D(X) ) = 1.  \eqno(1.0)
$$

 \  {\bf We will formulate our condition only in the simple and so-called {\it natural, or equally constructive terms,} which may be generated
through the trajectory of the considered random field (r.f,), in contradiction to the  many previous works; and we find also as a consequence the
 sufficient conditions for the weak compactness for the sequence of r.f. in these spaces.  } \par
 \  {\bf   We investigate also  as a capacity of an application of obtained results the classical Central Limit Theorem (CLT) in these spaces. } \par

\vspace{4mm}

 \   As for the previous works: see for example [2],[3],[4],[5],[6],[7],[10],[19],[20],
   [25],[27],[30], [36],[37]-[38],[42],[43] etc. \par

 \ The natural approach for the investigation of {\it continuous} random fields may be found in the articles [39],[12],[14],[26],[28]-[35] and so one.\par

\vspace{4mm}

 \ The well-known application of the CLT in the  Prokhorov-Skorokhod  spaces in the non-parametric statistics were obtained in the classical works of \\
Yu.V.Prokhorov  \ [40],  A.V.Skorokhod \  [43], see also  [25]. The multivariate generalization is considered in the
articles [3],[5],[7],[27],[37]-[38],[42] and so one. The important
application in the multi-parametric Monte-Carlo method may be found in [18],[20]. A very interest application of
 these limit theorems in physics  was investigated in [11],[37],[38]. \par
 \ The immediate predecessor of offered report is the preprint [30], in which was obtained the exact bilateral exponential bound for the
used in practice the tail probability

$$
{\bf P}_{\xi,X}(u) \stackrel{def}{=}  {\bf P} \left(\max_{x \in X} |\xi(x)| > u  \right), \ u > u_0 = \const \ge 1,
$$
for {\it discontinuous} random field $ \xi = \xi(x), \ x \in X. $  The one-dimensional case $ \ d = 1 \ $ was investigated in [36]. \par

\vspace{4mm}

  \ The paper is organized as follows. The section having even numbers: second, fourth  and sixth contains used further auxiliary apparatus.
 The third section is devoted  to the classical Lebesgue-Riesz approach for considered in this preprint problem. We offer in the fifth section
a more general method based on the theory  of the so-called Grand Lebesgue Spaces. \par

 \ The investigation of the Central Limit Theorem in the Prokhorov-Skorokhod spaces in the natural terms is the subject of the seventh
paragraph. We represent in the last section some concluding remarks. \par

\vspace{4mm}

 \ We must introduce some notations and definitions. Let $   \ x = \vec{x} = \{ x_j \} $ and $ \ y = \vec{y}   $ be two vectors from the source space $ \ X. \ $
Let also  $  \  q = q(x,y) \   $ be certain non zero  non - negative numerical values {\it continuous} symmetric function such that $  q(x,x) = 0, \ x \in X, $ and
let  $   x^+ = \vec{x}^+ \ge x, \ x^- = \vec{x}^- \le x,  $ where the inequalities are understood coordinate-wise:

$$
x,y \in X, \   x = \{ \ x(i) \ \}, \ y = \{ \ y(i) \ \}, \ x < y \ \Leftrightarrow \forall i = 1,2, \ldots,d \  x(i) < y(i).
$$
 \ Analogously may be defined the inequality  $ \ x \le y  $ and so one. \par

 \ We do not suppose that the function $ d(\cdot,\cdot) $
satisfies the triangle inequality. For instance, the function $ q(x,y)  $ may has a form

$$
q(x,y) = C | x - y|_{\alpha} \stackrel{def}{=}  C|x - y|^{\alpha}, \ \alpha = \const > 0,
$$
where $ \ |x| \ $ denotes the ordinary Euclidean norm of the vector $  x. $  Other example:

$$
C^{-1} \ q(x,y) =  \sum_{j=1}^d |x_j - y_j|^{\alpha(j)}, \ \alpha(j) = \const > 0.
$$

 \ We will name such a function  $  q = q(x,y) $ as a {\it quasy-distance.} \par
 \ We can and will suppose in the sequel without loss of generality $  \max_{x,y \in X} q(x,y) = 1.  $ \par

 \  These types spaces are introduced, applied and investigated in a recent report of Daniel J.Greenhoe [21], where was named as
"distance spaces".  \par

 \ The "metric" covering numbers $ \  N(X,q, \epsilon), \ \epsilon \in (0,1)  \  $  of the set $ \ X \   $ relative this  quasy-distance
 $  q = q(x,y) $ is defined  quite analogously to  classical metric-distance case, namely, as a minimal numbers of closed $ \  q -  \ $ "balls"

$$
B_q(x,\epsilon) := \{ y, \ y \in X, \ q(x,y) \le \epsilon  \}
$$
which cover all the set $  X: $

$$
 N(X,q, \epsilon) = \inf \{N: \ \exists x_k \in X,  k = 1,2,\ldots,N: \ \cup_{k=1}^N B_q(x_k, \epsilon) = X \}.
$$

 \ The natural logarithm of  the quantity $ \ H(X,q,\epsilon) := \ln  N(X,q, \epsilon)   \ $ is named as metric entropy of the set $  X  $ relative the
quasy-distance  $ q $ at the point  $ \epsilon, \ \epsilon \in (0,1). $ \par
 \ If for instance,

$$
q(x,y) = q_{\alpha}(x,y)  \asymp | x - y|^{\alpha}, \ \alpha = \const \in (0, \infty),
$$
 then

$$
N \left( [0,1]^d,  q_{\alpha}, \epsilon \right) \asymp \epsilon^{-d/\alpha}, \ \epsilon \in (0,1].
$$
 \ This notion play a very important role in the investigation of continuous random fields, see [12],[13]-[15],[28],[44]-[45].\par

\vspace{4mm}

 \ Further, let $  x^{(1)} = \vec{ x}^{(1)}, \   x^{(2)} = \vec{ x}^{(2)}, \  x^{(3)} = \vec{ x}^{(3)} \   $ be three (deterministic)
vectors from the set $  X  $ such that  $  \ x^{(1)}   \le  x^{(2)}  \le   x^{(3)}. \  $  \  Denote by $  T  $ the set of all the indexes
$  T = \{  1,2, \ldots, d   \} $  and by $  M  $  arbitrary subset of $  T:  $

$$
M = \{  i(1), i(2),  \ldots, i(m)   \}, \ 1 \le i(1) < i(2) < i(m) \le d. \eqno(1.1)
$$
 \ Of  course, $  m = 0 \ \Rightarrow M = \emptyset $ and $   m = d \   \Rightarrow M = T. \ $ Define the following vector $  z_M = \vec{z}_M  =
\vec{z}_M \left(x^{(1)}, \  x^{(2)} , \  x^{(3)} \right) \   $   generated by means of the random vectors $ x^{(1)}, \  x^{(2)} , \  x^{(3)} \ $  and by the
subset set $  M: $

$$
z_M(j) = x_j^{(1)}, \ j \notin M;  \ z_M(j) = x^{(3)}_j, \ j \in M. \eqno(1.2)
$$
 \ Evidently, $ z_X = x^{(1)}  $ and $  z_{\emptyset} = x^{(3)}.  \ $    Notice that the vector $  z_M $ dependent on the whole triple
$ \ \left(x^{(1)}, \  x^{(2)} , \  x^{(3)} \right). \   $ \par

 \ Define also for the source random field $ \xi = \xi(x) $ and for the certain ordered  triple $ \  \left(x^{(1)}, \  x^{(2)} , \  x^{(3)} \right) \   $
the system of partial increments

$$
\Delta[\xi](M) = \Delta[\xi](M)  \left(x^{(2)} \right)  :=  \left\{ \  \left( \xi(x^{(2)}) - \xi(z_M) \ \right ) \ \right\}, \ M \subset(T),
$$

$$
\tau[\xi] = \tau \left(x^{(1)}, \  x^{(2)} , \  x^{(3)} \right)[\xi] =  \min_{M \subset T} \left| \Delta[\xi](M)] \right|  =
$$

$$
 \min_{M \subset T} \left| \xi(x^{(2)}) - \xi(z_M)  \right|.  \eqno(1.3)
$$

 \ The generalized Prokhorov-Skorokhod   $   \kappa[\xi](h) =  \  \kappa[\xi]_q(h), \ h \in [ 0,1] \  $ module of continuity for the  random
field $  \xi(\cdot) $  may be defined as follows.   $ \  \kappa[\xi](h) = $

$$
\kappa[\xi]_q(h) \stackrel{def}{=} \sup_{x^{(2)} \in X} \sup_{q(x^{(3)}, x^{(1)}) \le h} \tau \left(x^{(1)}, \  x^{(2)} , \  x^{(3)} \right)[\xi] . \eqno(1.4)
$$

 \vspace{4mm}

 \ By definition, the random field $  \xi = \xi(x) $ belongs to the (multidimensional, in general case) Prokhorov-Skorokhod space $  \ D[0,1]^d \ $ almost
everywhere, iff for some (equally, each) non zero quasy-distance function $  \ q = q(x,y), \ x,y \in [0,1]^d $

$$
{\bf P} \left( \lim_{h \to 0+} \kappa[\xi]_q(h) = 0  \right) = 1. \eqno(1.5)
$$

\vspace{4mm}

 \ The complete investigation of these spaces in the multivariate case $  d \ge 2, $  in particular,  the  {\it criterion} for the tightness for the Borelian
measures in these spaces, may be found in [3],[5],[10],[27],[36],[37],[42] etc.\par

\vspace{4mm}

\section{Auxiliary estimates I. General approach.}

 \vspace{3mm}

 {\bf A.}  \  We will denote as customary for arbitrary r.v. $  \eta $  its  classical  Lebesgue - Riesz $  \ L_p = L_p(\Omega) \ $ norm

$$
|\eta|_p = |\eta|L_p:= \left[ {\bf E} |\eta|^p \right]^{1/p}, \ p = \const \ge 1.
$$

 \  A multivariate generalization of this notion may be introduced as follows.
  \ Let  $  \xi = \vec{\xi} = \{ \xi_1, \xi_2,\ldots, \xi_k  \}  $  be a random vector and let $  \ p = \vec{p} = \{  p_1, p_2,\ldots, p_k  \} \ $ be
numerical deterministic $  \ k \ - \ $ dimensional vector  such that $  \ \forall i \ \Rightarrow p_i \ge 1. $
 We introduce  the following  {\it mixed}  moment

$$
\mu = \mu \left(  \{  \xi_i \}, \{p_i \}  \right) =  \mu \left( \ \vec{\xi}, \  \vec{p} \ \right) =
$$

$$
\mu[\xi](p_1, p_2, \ldots, p_k ) \stackrel{def}{=}  {\bf E} \left[ \prod_{i=1}^k  \left|\xi_i \right|^{p_i}  \right].
$$

 \ Let $  \{ \ a(1), a(2), \ldots, a(k)  \ \}  $ be again $ \ k - \ $ tuple of real numbers  greatest that one:  $  a(j) > 1, $
and such that

$$
\sum_{i=1}^k \frac{1}{a(i)} = 1.
$$
 \ The set all of such the $ \ k \ $  tuples we will denote by $  A = A(k). $   We apply  the famous H\"older's inequality

$$
\mu \le \prod_{j=1}^k  \  \left[ {\bf E} |\xi_i|^{a(i) \cdot p_i} \right]^{1/a(i)}  =
 \prod_{i=1}^k \left\{ \left| \xi_i \ \right|_{a(i) \ p_i} \right\}^{p_i}. \eqno(2.A)
$$

  \ The last estimate may can be strengthened as follows.

$$
\mu  \le \inf_{\vec{a} \in A(k)} \left[ \ \prod_{j=1}^k  \  \left[ {\bf E} |\xi_i|^{a(i) \cdot p_i} \right]^{1/a(i)} \ \right]  =
 \inf_{\vec{a} \in A(k)} \left[ \left[ \ \prod_{i=1}^k    \left| \xi_i \ \right|_{a(i) \ p_i} \right\}^{p_i} \ \right]. \eqno(2.B)
$$

\vspace{4mm}

{\bf B.} \ Denote $ \ \zeta = \ $

$$
 \zeta_q [\xi] \left(x^{(1)}, \ x^{(3)}; \ u \right)  = \zeta [\xi] \left(x^{(1)}, \ x^{(3)}; \ u \right) :=
\sup_{x^{(2)} \in X}  {\bf P} \left( \ \min_{M \subset T}  |\Delta[\xi](M) | > u \ \right) =
$$

$$
 \sup_{x^{(2)} \in X}  {\bf P} \left(  \min_{M \subset T}  \left|  \ \xi(x^{(2)}  - \xi(x_M)  \  \right| \  > u  \right),
$$

$$
x^{(1)}, \ x^{(2)}, \  x^{(3)} \in X, \hspace{4mm} u > 0.
$$

 \ We will start from the following condition imposed  by the  natural way on the paths of the random field $ \xi(\cdot). $

\vspace{3mm}

 \ {\bf Condition 2.0.} Assume that for certain continuous quasy-distance $  q = q(x,y), \ x,y \in X  $

$$
  \zeta_q[\xi] \left(x^{(1)}, \ x^{(3)}; \ u \right) \le \frac{q(x^{(1)}, x^{(3)})}{\lambda(u)}, \ u \ge 1, \eqno(2.0)
$$
 is the so-called condition of {\it key estimate},  where $ \lambda = \lambda(u) $ is certain monotonically increasing function,
continuous or not, satisfying the  condition $ \ \lim_{u \to \infty} \lambda(u) = \infty. $ \par

\vspace{3mm}

{\bf Remark 2.1.} \ Both the functions  $ q(\cdot, \cdot) $ and $ \lambda(\cdot) $  in (2.0) may be introduced by the natural way as follows.

$$
q(x^{(1)}, x^{(3)} ) = q_{\xi}(x^{(1)}, x^{(3)})
\stackrel{def}{=} \sup_{x^{(2)} \in X} {\bf P} \left[ \ \tau \left(x^{(1)}, \  x^{(2)} , \  x^{(3)} \right) \ge 1 \ \right] \eqno(2.1)
$$
and

$$
\frac{1}{\lambda(u)} = \frac{1}{\lambda_{\xi}(u)} \stackrel{def}{=}
\sup_{0 \le x^{(1)} < x^{(3)}} \ \left[ \frac{ {\bf P} \left[ \ \tau \left(x^{(1)}, \  x^{(2)} , \  x^{(3)} \right) \ge u \ \right] }{ q(x^{(1)}, x^{(3)} )  } \right].
\eqno(2.2)
$$

 \ Herewith the estimate (2.0) there holds, if obviously both these functions $  \ q(\cdot, \cdot) , \ \lambda(\cdot) $ there exist. \par

\vspace{4mm}

 \ {\bf C.} \  An alternative but again constructive method for quasy-distance estimation (2.0) based of the classical Lebesgue-Riesz spaces
$ \ L_p = \ L_p(\Omega) \ $ is  follows. Introduce the following vectors and functions:

$$
s = \vec{s} = \{ \ s(M)  \  \}, \ M \subset T,  s(M)  > 0.
$$
 \ The set of such a vectors  $ \   \{s = s(M) \} \  $ will be denoted by $  L; \ L = \{  \vec{s}\}. $ \par

\ Further, introduce the following function

$$
 \beta \left( \ x^{(1)}, \ x^{(3)} \ \right) =  \beta \left(x^{(1)}, \ x^{(3)}; \vec{s} \right)  =
 \beta[\xi] \left(x^{(1)}, \ x^{(3)}; \vec{s} \right) \stackrel{def}{=}
$$

$$
\sup_{x^{(2)} \in X} {\bf E} \prod_{M \subset T} \left[ \  \left|  \ \Delta[\xi](M) \right|^{s(M)} \ \right]. \eqno(2.3)
$$

 \ The function $ \ \beta(x,y), \ x,y \in X \ $ by fixed value of the vector $  \ \vec{s} \ $ is quasy-distance function
on the set $  X \otimes X,  $  if of course  there exists for some positive  vector $  \ \vec{s}. $\par

 \ We derive using the famous Tchebychev's inequality the restriction of type (2.0)

$$
 {\bf P}_{\zeta[\xi]}(x^{(1)}, \  x^{(3)}, \ u) \le \inf_{x^{(2)} \in X}  \frac{\beta \left(x^{(2)}, \vec{s} \right)}{u^{ \sum_M s(M)  }}, \eqno(2.4)
$$
if  evidently the right-hand side of the last inequality is finite. So, in this case $  \lambda(u) = u^{\sum_M s(M)}, \ u > 0. $  \par

 \  The right-hand side of the inequality (2.4) may be estimated in particular as follows.

$$
 {\bf P}_{\zeta[\xi]}(x^{(1)}, \  x^{(3)}, \ u) \le \sup_{x^{(2)} \in X}  \left[ \ \frac{\beta \left(x^{(2)},  t \cdot \vec{1} \right)}{u^{ 2 d t }} \ \right], \eqno(2.4a)
$$
$    \vec{1} :=  \{ 1,1, \ldots, 1   \};  $  \ \hspace{5mm}  $   t = \const > 0. $\par

 \ Obviously,

$$
 {\bf P}_{\zeta[\xi]}(x^{(1)}, \  x^{(3)}, \ u) \le  \inf_{t > 0} \ \sup_{x^{(2)} \in X} \ \left[ \frac{\beta \left(x^{(2)},  t \cdot \vec{1} \right)}{u^{ 2 d t }}\right]. \eqno(2.4b)
$$

 \ Se for the detail explanation [36].\\

 \ Somehow or other, we grounded the key estimate (2.0). \par

 \ Alike estimate with replacing the classical Lebesgue-Riesz space $  L_p $ by so - called Grand Lebesgue Space will be considered further. \par

\vspace{4mm}

 \ Denote by $  \varepsilon $ the set of all positive sequences $  \   \varepsilon = \{\epsilon(k) \} $ such that $  \epsilon(1) = 1, \  \epsilon(k) \downarrow 0, \ $ and
by $ \Theta  $ the set of all positive sequences $  \Theta = \{   \theta(k)  \} $ for which $ \theta(k) \downarrow 0 $ and

$$
\sum_{k=0}^{\infty} \theta(k) = 1.
$$

 \vspace{3mm}

 \ Introduce the following variables:

$$
Q = Q(X,q, \xi;u) \stackrel{def}{=} \inf_{ \{\epsilon\} \in \varepsilon  } \inf_{ \{ \theta \} \in \Theta  }
\sum_{k=0}^{\infty} N ( X,q,\epsilon(k+1) ) \cdot \frac{\epsilon(k)}{\lambda(u \cdot \theta(k))}.
$$

$$
\sigma[q](h) \stackrel{def}{=} h^{-d} \sup_{ |x^{(3)} - x^{(1)}| \le 2 h} q \left(x^{(3)}, x^{(1)} \right).
$$

 \ We need to use now  a very important fact, which is an evident and slight generalization of the main result of the articles [30],[36]; see also the predecessor [42].\\

\vspace{4mm}

{\bf Theorem 2.1.}   \\

{\bf A.} Suppose in addition to the assumption (2.0) that  for the separable stochastic continuous random field $ \ \xi = \xi(x), \ x \in X =[0,1]^d $ the
following conditions holds true:

\vspace{3mm}

$$
\lim_{u \to \infty}  Q(X, q, \xi; u) =0;
$$
then

\vspace{3mm}

$$
{\bf P}(\tau[\xi] > u) \le  Q(X,q, \xi;u). \eqno(2.5)
$$

\vspace{3mm}

 {\bf B.}  As long as

\vspace{3mm}

$$
\lim_{h \to 0+} \sigma[q](h)  = 0,
$$
we deduce

\vspace{3mm}

$$
{\bf P}(\kappa[\xi](h) > u) \le   Q(X,q, \xi;u) \cdot \sigma[q](2 h), \eqno(2.6)
$$
and as a consequence  $  {\bf P} \left(\xi(\cdot) \in D[0,1]^d = 1 \right).  $\par

\vspace{4mm}

 \ {\bf Example 2.1.}  Suppose

$$
\lambda(u) = u^{2 \rho}, \  u \ge 1, \ \exists  \rho = \const  > 0,
$$
and that the quasy - distance $  q(\cdot, \cdot)  $ is such that

$$
N(X, q, \epsilon) \le C_N \ \epsilon^{-\gamma}, \  \gamma = \const \in (0,1).
$$
 \ One can apply the statement of theorem 2.1 choosing correspondingly

$$
\epsilon(k) := s^{k-1}, \  \theta(k) = (1 - s) s^k, \ s = \const \in (0,1):
$$

$$
 \frac{Q(X,q, \xi;u)}{ C_N \ u^{-2 \rho}  } \le \frac{s^{-\gamma} \ (1 - s)^{-2 \rho}}{1 - s^{ 1 - s^{1 - \gamma - 2 \rho} }}, \
\gamma + 2 \rho < 1,
$$
therefore

$$
 \frac{Q(X,q, \xi;u)}{ C_N \ u^{-2 \rho}  } \le \inf_{s \in (0,1)}
\left[ \ \frac{s^{-\gamma} \ (1 - s)^{-2 \rho}}{1 - s^{ 1 - s^{1 - \gamma - 2 \rho} }} \ \right],
$$
so that

$$
Q(X,q, \xi;u) \le     C(C_N, \gamma, \rho) \  u^{-2 \rho},  \ u \ge 1, \eqno(2.7)
$$
if of course $ \ \gamma + 2 \rho < 1. $ \par

 \  By virtue of (uniform) continuity of the  quasy-distance function $  \ q(\cdot, \cdot) $

$$
\lim_{h \to 0+} \sigma[q](h)  = 0,
$$
following

$$
{\bf P}(\kappa[\xi](h) > u) \le C(C_N, \gamma, \rho) \  u^{-2 \rho}  \cdot \sigma[q](2 h), \ u \ge 1, \eqno(2.8)
$$
and as a consequence  $  {\bf P} \left(\xi(\cdot) \in D[0,1]^d \right) = 1.  $\par

\vspace{4mm}

 \  Assume  in continuation that the {\it family,} or for simplicity the {\it sequence} of somehow dependent random fields
$  \xi_n(X), \ n = 1,2,\ldots $ be a given such that all the introduced conditions are satisfied uniformly on the parameter $ n. $  \par

 \ In detail, suppose

$$
\sup_n {\bf P}_{\zeta [\xi_n]}(u) = {\bf P}_{\kappa[\xi_n]}(x^{(1)}, \  x^{(3)}; u) \stackrel{def}{=}
$$

$$
\sup_n \sup_{x^{(2)}  \in X}  {\bf P} \left[ \ \tau[\xi_n] \left(x^{(1)}, \  x^{(2)} , \  x^{(3)} \right) \ge u \ \right] \le
$$

$$
\frac{q(x^{(1)}, x^{(3)} )}{\lambda(u)}, \ u \ge 1, \eqno(2.9)
$$
which may be called as the {\it uniform key estimate},  where $ \lambda = \lambda(u) $ is a monotonically decreasing function, continuous or not,
with condition $ \ \lim_{u \to \infty} \lambda(u) = \infty. $ \par

\vspace{4mm}

{\bf Theorem 2.2.}   Let the finite - dimensional distribution of the r.f. $ \xi_n(\cdot) $ converges to ones for some r.f. $  \xi_{\infty}(\cdot) $
belonging to at the same space $  \ D[0,1]^d \ $ with probability one. \par

 \ Suppose that  for the random fields $ \ \xi_n = \xi_n(x), \ x \in X =[0,1]^d $ the following uniform
conditions holds true:

$$
\lim_{u \to \infty} \sup_n  Q(X,q, \xi_n ;u) =0;  \eqno(2.10)
$$
then

$$
 \sup_n {\bf P}(\Delta[\xi_n] > u) \le  \sup_n Q(X,q, \xi_n; \ u) =: \overline{Q} (X,q, \{\xi_n(\cdot);\} \ u). \eqno(2.11)
$$

  \ Since

$$
\lim_{h \to 0+} \sigma[q](h)  = 0, \eqno(2.12)
$$
we obtain

$$
 \sup_n {\bf P}(\kappa[\xi_n](h) > u) \le  \sup_n Q(X,q, \xi;u) \cdot \sigma[q](2 h), \eqno(2.13)
$$
and as a consequence  $  {\bf P} \left(\xi_n(\cdot) \in D[0,1]^d = 1 \right)  $ and furthermore the sequence of the  r.f. $  \xi_n(\cdot) $ converges
in distribution  in the space   $ D[0,1]^d $ to the distribution of the r.f. $  \xi_{\infty}. $\par

\vspace{4mm}

{\bf The Proof} follows immediately from the estimate

$$
\sup_n {\bf P}(\kappa[\xi_n](h) > u) \le  \sup_n  \overline{Q}(X,q, \{ \xi_n\}; \ u) \cdot \sigma[q](2 h),
$$
therefore the sequence of distributions of the sequence $  \ \xi_n(\cdot) $ is weakly compact in the considered space.
See in detail  [3], [9], [10], [27], [42]. \par

 \ Evidently, instead the integer index $  \ n \ $ may be used a point  $ \alpha  $ belonging on arbitrary set, in particular, some topological space. \par

 \vspace{3mm}
\section{Classical Lebesgue-Riesz approach.}
 \vspace{3mm}

 \ We rely here on the  introduced before the generalized Lebesgue-Riesz distance  $  \ \beta(x,y) = \beta(x,y, \vec{s}),  \ x,y \in X = [0,1]^d $  and
the statement of theorem 2.1. \par

 \ We adopt in this section the vector  $ \ s =  \vec{s} $ to be  {\it arbitrary} positive fixes constant vector.   \par

 \ Assumptions and notations:

$$
  p := \sum_M s(M) > 0, \ \gamma = \const  \in (0,1), \ C(N)  = \const \in (0,\infty); \ \forall \epsilon \in (0,1) \Rightarrow
$$

$$
N(X,\beta(\cdot, \cdot), \ \epsilon) \le C(N) \ \epsilon^{-\gamma},   \eqno(3.0)
$$
following,
$$
\lambda(u) = C(\lambda) \  u^{\sum_M s(M) }= C(\lambda) \  u^p ,  \ u > 0;  \ C(\lambda) = \const \in (0,\infty).
$$

 \ One can choose in the statement of theorem 2.1

$$
\epsilon(k) = 2^{-k}, \hspace{5mm} \theta(k) = (1 - \theta_0) \ \theta^k_0,
$$
where $   \theta_0 = 2^{  (\gamma - 1)/2p } \in (0,1). \ $  Denote

$$
 W(\gamma) := \left( 1 - 2^{  - (1 - \gamma)/2}  \right)^{-1}.
$$

\vspace{4mm}

{\bf Theorem 3.1.} We get under our notations and conditions after some computations:

 $$
  Q(X,\beta,\xi; \ u) \le \frac{2 C(N)}{C(\lambda)} \ W^{p + 1}(\gamma) \ u^{-p}, \ u \ge 1. \eqno(3.1)
 $$

 \ If in addition

$$
\lim_{h \to 0+} \sigma[\beta](h)  = 0, \eqno(3.2)
$$
then

$$
{\bf P}(\kappa[\xi](h) > u) \le \ \frac{2 C(N)}{C(\lambda)} \ W^{p + 1}(\gamma) \ u^{-p}  \cdot \sigma[\beta](2 h), \ u \ge 1, \eqno(3.2)
$$
and as a consequence  $  {\bf P} \left(\xi(\cdot) \in D[0,1]^d \right) = 1.  $\par

\vspace{4mm}

 \ {\bf Example 3.1.} The conditions of theorem (3.1)  are trivially satisfied for the r.f. $ \  \xi =\xi(x), x \in X \ $ of the form

$$
\xi(x) = I(\vec{\eta} < \vec{x}), \eqno(3.3)
$$
where $  \ \vec{\eta} \  $ is a random vector with values in the set $ \  X  \ $ having  at last continuous function of distribution

$$
F(x) = F(\vec{x}) = {\bf P} (\vec{\eta} < \vec{x}); \eqno(3.4)
$$
the inequality between two vectors is understood as ordinary coordinate-wise, and the $ \ I(\cdot) \  $ is ordinary indicator function.\par

 \ Analogous statement is true under at the same assumptions  for the r.f.

$$
\xi(x) = I(\vec{\eta} < \vec{x})  - F(\vec{x}). \eqno(3.5)
$$

 \ The grounding follows immediately from the relation

$$
\beta \left(x^{(1)},  x^{(2)},  x^{(3)} ; \vec{s} \right) = 0,
$$
if $ \ \min_j x^{(1)}_j > 0 \ $ or  $ \ \max_j x^{(3)}_j < 1 . $ \par

\vspace{3mm}

 \ {\bf Remark 3.1.} Note that  the variable $ \ \beta(\cdot, \cdot; \cdot) \ $ allows by virtue of H\"older's inequality the following simple estimate

$$
\beta \left(x^{(1)}, x^{(3)}; \vec{s} \right) \le
$$
$$
\sup_{x^{(2)} \in X} \prod_{M \subset T} \left[ {\bf E} | \Delta[\xi](M)  |^{\alpha(M) \ s(M)} \right]^{1/s(M)} =
$$

$$
\sup_{x^{(2)} \in X}  \prod_{M \subset T} |  \Delta[\xi](M) |^{s(M)}_{\alpha(M) \ s(M)}, \eqno(3.6)
$$
where
$$
\alpha(M) > 1, \ \sum_M 1/\alpha(M) = 1.
$$

 \vspace{3mm}

\section{ Auxiliary facts II: Grand Lebesgue Spaces (GLS).}

 \vspace{3mm}

 \ We recall here some facts about the so-called  Grand Lebesgue  Spaces (GLS) spaces and deduce also mixed norm inequalities for the random
vector belonging to these spaces. \par

\vspace{3mm}

 \ Let  $  \psi = \psi(p), \ p \in [1, b), \ b = \const \in (1,\infty]  $ (or   $ p \in [1,b] $ ) be certain bounded
from below:  $  \ \inf \psi(p)  > 0 $ continuous inside the  {\it  semi - open} interval $   \ p \in [1, b) $ numerical function. We can and will suppose
$   \ b = \sup{p, \psi(p) < \infty}, $ so that  $    \supp \psi = [1, b)  $  or $ \supp \psi = [1, b]. $ The set of all such a functions will be denoted by
$ \Psi(b) = \{  \psi(\cdot)  \}; \ \Psi := \Psi(\infty).  $\par

 \ By definition, the (Banach) Grand Lebesgue Space  \ (GLS)  \ space   $  \ G\psi = G\psi(b)  $ consists on all the numerical
valued random variables   $ \zeta  $ defined on our probability space  and having a finite norm

$$
||\zeta|| = ||\zeta||G\psi \stackrel{def}{=} \sup_{p \in [1,b)} \left\{ \frac{|\zeta|_p}{\psi(p)} \right\}. \eqno(4.0)
$$

 \ These spaces  are Banach functional space, are complete, and rearrangement invariant in the classical sense, see [1], chapters 1,2;
and were  investigated in particular in  many  works, see e.g.[16]-[17],[23],[26],[28],[35],[36].
 We refer here  some  used in the sequel facts about these spaces and supplement more. \par

\vspace{4mm}

{\bf Remark 4.1.}  \ Let us consider the so-called {\it degenerate} $  \ \Psi \ - \ $ function $ \ \psi_{(r)}(p), \ $ where $  r = \const \in [1,\infty): $

$$
\psi_{(r)}(p) \stackrel{def}{=} 1,  \ p \in [1,r];
$$
so that the correspondent value $  b = b(r) $  is equal to $  r. $  One can  extrapolate formally this function onto the whole semi-axis $  R^1_+: $

$$
\psi_{(r)}(p)  := \infty, \ p > r.
$$

 \  The classical Lebesgue-Riesz $ L_r  $ norm for the r.v. $  \eta $ is  quite equal to the GLS norm $  ||\eta|| G\psi_{(r)}: $

$$
|\eta|_r = ||\eta|| G\psi_{(r)}.  \eqno(4.1)
$$
 \ Thus, the ordinary Lebesgue-Riesz spaces are particular, more precisely, extremal cases of the Grand-Lebesgue ones. \par

\vspace{4mm}

\  Suppose now  $  0 < ||\zeta || < \infty. $
Define the function $  v(p) = v_{\psi}(p) := p  \ \ln \psi( p), \ 1 \le p < b $ and put formally $ ν(p) := \infty, \ p < 1 $ or $ p > b. $  \\

 \ Recall that the Young-Fenchel,or Legendre transform $  f^*(y) $ for arbitrary function $ f : R \to R $ is defined (in the one-dimensional case)
as follows

$$
 f^*(y)  \stackrel{def}{=} \ \sup_{x \in \Dom(f)} (xy - f(x)).
$$

 \  It is known that

$$
{\bf P}(|\zeta| > y) \le \exp \left(-v_{\psi}^* (\ln (y/||\zeta||) \right), \ y \ge e \cdot ||\zeta||. \eqno(4.2)
$$

  \ Conversely, the last inequality may be reversed  in the following version: if

$$
{\bf P}(|\zeta| > y) \le \exp \left(-v_{\psi}^* (\ln (y/K) \right). \ y \ge e \cdot K, \ K = \const \in (0,\infty),
$$
and if the function $ v_{\psi}(p), \ 1 \le p < \infty  \ $  is positive, continuous, convex and such that

$$
\lim_{p \to \infty} \ln \psi(p) = \infty,
$$
then  $ \zeta \in G\psi $ and besides $  ||\zeta||  \le C(\psi) \cdot K.  $\par

 \ Moreover, let us introduce the so-called  {\it exponential} Orlicz space $ \ L(M ) \ $  over the source probability space with
the generating Young-Orlicz function

$$
M(u) := exp \left( v_{\psi}(\ln |u|) \right), \ |u| > e
$$
and as ordinary $ M(u) = \exp(C u^2) - 1, \ |u| \le e. $ It is known, [26],[35],[36] that the
$  G\psi $ norm of arbitrary r.v. $ \zeta $ is  equivalent to the its norm in the (exponential) Orlicz space  $ L(M ): $

$$
||\zeta||G\psi \le C_1 ||\zeta||L(M)  \le C_2 ||\zeta||G\psi, \ 0 < C_1 < C_2 < \infty. \eqno(4.3)
$$

 \  Furthermore,  let now $  \eta = \eta(z), \ z \in Z $ be arbitrary family  of random variables  defined on any set $  Z  $ such that

$$
\exists b \in (1,\infty] \ \forall p \in [1,b)  \ \Rightarrow  \psi_Z(p) := \sup_{z \in Z} |\eta(z)|_p  < \infty. \eqno(4.4)
$$
 \ The function $  p \to \psi_Z(p)  $ is named as {\it  natural} function for the  family  of random variables $  Z.  $  Obviously,

$$
\sup_{z \in Z} ||\eta(z)||G\Psi_Z = 1.
$$

\vspace{3mm}

  \ Let again  $  \xi = \vec{\xi} = \{ \xi_1, \xi_2,\ldots, \xi_k  \}  $  be a random vector and let $  \ p = \vec{p} = \{  p_1,p_2,\ldots, p_k  \} \ $ be
numerical deterministic $  \ k \ - \ $ dimensional vector  such that $  \ \forall i \ \Rightarrow p_i \ge 1. $
 We introduced in the second section  the following  {\it mixed}  moment

$$
\mu = \mu \left(  \{  \xi_i \}, \{p_i \}  \right) =  {\bf E} \left[ \prod_{i=1}^k  \left|\xi_i \right|^{p_i}  \right]
$$
and obtained the following estimate

$$
\mu  \le  \ \inf_{\vec{a} \in A(k)}  \left[ \ \prod_{i=1}^k  \left[  \left| \xi_i \ \right|_{a(i) \ p_i} \right]^{p_i} \ \right].
$$

 \ Suppose now that each r.v.  $   \xi_i $  belongs to certain $ \ G\psi_i =  G\psi_i(b(i)), \ 1 < b(i) \le \infty  $ space:

$$
|\xi_i|_p \le  ||\xi_i||G\psi_i \cdot \psi_i(p), \ 1 \le p \le \infty;
$$
then

$$
 \ \mu \le Y = Y(  \{  ||\xi_i||G\psi_i \},   \{ p_i  \}   ),   \eqno(4.5)
$$
where

$$
 Y = Y(  \{  ||\xi_i||G\psi_i \},   \{ p_i  \} ) \stackrel{def}{=}
$$
$$
  \ \inf_{\vec{a} \in A(k)}  \ \prod_{i=1}^k \left[ \    ||\xi_i||G\psi_i \cdot  \psi_i \left( a(i) \ p_i \right) \ \right]^{p_i}.  \eqno(4.6)
$$

 \  If for instance

$$
\forall i \ ||\xi_i||G\psi_i =1,
$$
for example when all  the functions $  \psi_i $ are the natural function for the r.v. $ \ \xi_i, \  $ then the estimate (4.6)  may be evidently simplified
as follows

$$
\ \mu \le \ \inf_{\vec{a} \in A(k)}  \ \prod_{i=1}^k \left[ \   \left( a(i) \ p_i \right) \ \right]^{p_i}.  \eqno(4.6a)
$$

\ We deduce the following estimate for  the tail of minimum of $ \ \xi_i: $

$$
{\bf P} \left( \min_i |\eta_i| > u  \right) \le \frac{Y(  \{  ||\xi_i||G\psi_i \},   \{ p_i  \}   )}{u^{\sum_i p_i}}. \eqno(4.7)
$$

 \ This relation play a very important role for the investigation of discontinuous random fields and will be used further. \par

 \vspace{3mm}

\section{Grand Lebesgue Spaces approach.}

\vspace{3mm}

 \ We rely here on the introduced before the generalized Lebesgue-Riesz distance  $  \ \beta(x,y) = \beta(x,y, \vec{s}),  \ x,y \in X = [0,1]^d $  and
the statement of theorem 2.1 to extend the obtained therein results into the so-called {\it \ Grand Lebesgue Spaces (GLS)}. \par

 \vspace{3mm}

 \ We adopt in the sequel the vector  $ \ s =  \vec{s} $ to be  {\it arbitrary, i.e. variable} positive vector.   \par

 \vspace{3mm}

 \ Assumptions:

$$
\exists \gamma(s) \in [0,\infty), \ C_s(N)  = \const \in (0,\infty), \ \forall  \epsilon \in (0,1) \ \Rightarrow
$$

$$
N(X,\beta(\cdot, \cdot; \ s), \ \epsilon) \le C_s(N) \ \epsilon^{-\gamma(s)};  \eqno(5.1)
$$
following,
$$
\lambda(u) = \lambda_s(u) = C_s(\lambda) \  u^{\sum_M s(M)   } ,  \ u > 0;  \ C_s(\lambda)  = \const \in (0,\infty).  \eqno(5.2)
$$

 \ The subset of the whole set  $ \ L \ $  of such a vectors  $ \  s =  \{ s(M) \} \  $    for which $ \gamma(s)   \in [0,1) $
will be denoted by $  L^o; \ L^o = \{ s =  \vec{s}:  \ \gamma(s) \in [0,1); \ L^o  \subset L \}. $  {\it We suppose in the sequel} $  \ L^o \ne \emptyset. $  \par

\vspace{4mm}

  \ Let us return now to the source problem concerning the Prokhorov-Skorokhod continuity of the r.f. $  \ \xi = \xi(x). $ \par

   \ We get under our notations and conditions  by virtue of theorem 3.1 for all the values $  \ s \in L^o $

 $$
  Q(X,\beta(\cdot, \cdot; \ s), \ \xi; \ u) \le \frac{2 C_s(N)}{C_s(\lambda)} \ W^{p(s) + 1}(\gamma(s)) \ u^{-p(s)}, \ u \ge 1. \eqno(5.3)
 $$

 \ If in addition

 $$
  \exists \ s \in L^o \ \Rightarrow
\lim_{h \to 0+}  \sigma[\beta(\cdot, \cdot;  \ s)](h)  = 0,  \eqno(5.3a)
$$
i.e. if the function $  \ (x,y) \to \beta(x,y; s)  $ is (uniform) continuous, then   $   \ \exists\ s \in L^o  \ \Rightarrow $

$$
{\bf P}(\kappa[\xi](h) > u) \le \ \frac{2 C_s(N)}{C_s(\lambda)} \ W^{p(s) + 1}(\gamma(s)) \ u^{-p(s)}  \cdot \sigma[\beta(\cdot, \cdot;  s)](2 h), \ u \ge 1.
 \eqno(5.4)
$$

 \  Denote by $  L^+ $ the set of all the vectors $  \ s = \vec{s} \ $  from the set $  \ L^o  \ $   for which the relation (5.3a) there holds.
{\it It is natural to suppose  also} $ \  L^+ \ \ne \emptyset. $\par

 \ The estimate (5.4) can be  obviously strengthened as follows $ {\bf P}(\kappa[\xi](h) > u) \le \  $

$$
\inf_{ s \in L^+  } \left\{  \frac{2 C_s(N)}{C_s(\lambda)} \ W^{p(s) + 1}(\gamma(s)) \ u^{-p(s)}
 \cdot \sigma[\beta(\cdot, \cdot;  s)](2 h) \ \right\}, \ u \ge 1.  \eqno(5.4a)
$$

\vspace{3mm}

 \ To summarize: \par

\vspace{4mm}

{\bf Theorem 5.1.}  We conclude under conditions of this section  $ \ {\bf P}(\kappa[\xi](h) > u) \le  $

$$
 \inf_{s \in L^+} \ \left[ \frac{2 C_s(N)}{C_s(\lambda)} \ W^{p(s) + 1}(\gamma(s)) \ u^{-p(s)}  \cdot \sigma[\beta(\cdot, \cdot;  s)](2 h)  \right],  \ u \ge 1,   \eqno(5.5)
$$
and as a consequence  $  {\bf P} \left(\xi(\cdot) \in D[0,1]^d \right) = 1.  $\par

\vspace{4mm}

 \ {\bf Example 5.1.}  Denote

$$
K(q) := \sup_{ s \in L^o, p(s) \le q } \left[  \frac{2 C_s(N)}{C_s(\lambda)} \ W^{p(s) + 1}(\gamma(s)) \right], \ 1 < q  < \infty,
$$
then

$$
{\bf P}(\zeta[\xi] > u) \le K(q) \ u^{-q}, \eqno(5.6)
$$
hence

$$
{\bf P}(\zeta[\xi] > u) \le \inf_{q > 0} \left[ \ K(q) \ u^{-q} \right]. \eqno(5.7)
$$

 \  If in addition

$$
 \exists m = \const > 0 \  \forall q > 1 \   K(q) \le C_1 q^{ q/m },
$$
then

$$
{\bf P} \left(\zeta[\xi](x^{(1)}, \  x^{(3)} )> u \right) \le  C_1 \ \exp \left\{- u^m/(me) \right\}, \ u \ge e.
$$

 \ Thus, we obtained an {\it exponential}  bound for tail of distribution of  the r.v. $  \zeta[\xi],   $ in accordance with the
general theory of Grand Lebesgue Spaces, see fourth section.\par
 \ More simple result may be obtained from the inequality (2.4a):

$$
 {\bf P}_{\zeta[\xi]}(x^{(1)}, \  x^{(3)}, \ u) \le \inf_{t > 0}
\left[ \ \sup_{x^{(2)} \in X}  \frac{\beta \left(x^{(2)},  t \cdot \vec{1} \right)}{u^{ 2 d t }} \ \right]. \eqno(5.8)
$$

\vspace{4mm}

 \ {\bf Remark 5.1.} \ Appearing in (1.3),  (1.4) the r.v.

$$
\tau[\xi] = \tau \left(x^{(1)}, \  x^{(2)} , \  x^{(3)} \right)[\xi] =  \min_{M \subset T} \left| \Delta[\xi](M)] \right|
$$
may be estimated by means of the general theory of GLS as follows.\par

 \ Suppose that each of the r.v. $  \ \Delta[\xi](M), \ M \subset T \ $  belongs uniformly over the space $  X  $ to some $ \ G\psi_M \ $ - space:

$$
\forall M \subset T \  \exists b(M) \in (1, \infty]  \ \forall p \in [1,b(M) \Rightarrow   \psi_M(p) :=
$$

$$
 \sup_{x^{(1)},  x^{(2)},  x^{(3)} \in X} \left| \ \Delta_M \ \right|_p < \infty.
$$

 \ We define formally as before the values $  \psi_M(p) $ in the case when $  b(M) < \infty $  for the values $  \ p > b(M) \ $
 as $ \infty, $ so that for all the values $  p \ge 1 $

$$
 \psi_M(p) = \sup_{x^{(1)},  x^{(2)},  x^{(3)} \in X} \left| \ \Delta_M \ \right|_p.  \eqno(5.9)
$$

 \  We deduce by means of inequality (4.7)

$$
{\bf P} (\tau[\xi] > u) = {\bf P} \left( \min_{M \subset T} \left| \Delta[\xi](M)] \right| > u \right) \le
$$

$$
\inf_{ \{ p_M \ \}} \frac{Y \left( \ \left\{ ||s_M||\psi_M  \right\},  \ \left\{  p_M   \right\} \  \right)}{ u^{ \sum_{M \subset T } p(M)} }. \eqno(5.10)
$$

\vspace{4mm}

 \ {\bf Another approach.} \par

\vspace{4mm}

 \ Let us define this time the correspondent natural  $  \ \psi \ -  $ function $ \ \upsilon = \upsilon(p) \ $ as follows

$$
\upsilon(p) := \sup_{ x^{(1)}, x^{(2)}, x^{(3)} \in X  } \left| \ \tau \left(x^{(1)}, \  x^{(2)} , \  x^{(3)} \right) \ \right|_p, \eqno(5.11)
$$
if naturally $ \ \upsilon(p) \ $ is finite in some non-trivial closed or semi-closed  interval $ p \in [1, b) $ or  $ p \in [1, b], \ b = \const \in (1, \infty]. $ \par

 \ Another equivalent version:

$$
\upsilon(p) := \sup_{ x^{(1)}, x^{(2)}, x^{(3)} \in X  } \left| \ \Delta[\xi]  \left(M \left(x^{(1)}, \  x^{(2)} , \  x^{(3)} \ \right) \right) \ \right|_p. \eqno(5.11a)
$$

 \ The correspondent natural distance function $ r(x,y) \ $  takes here a form

$$
r = r(\cdot, \cdot) = r_{\upsilon} \left((x^{(1)}, \  x^{(3)}  \right) =  r \left((x^{(1)}, \  x^{(3)}  \right) \stackrel{def}{=}
$$

$$
 \sup_{x^{(2)} \in X} \left| \left|  \left[ \ \tau \left(x^{(1)}, \  x^{(2)}, \  x^{(3)} \right)  \ \right] \ \right| \right|G\upsilon, \eqno(5.12)
$$
and analogously may be defined  as before the variable $ \ \sigma_r(h). $ \par

 \  We impose  at the same condition

$$
\forall \epsilon \in (0,1) \   N(X,r,\epsilon) \le C(N,r)  \epsilon^{-\gamma},  \ \exists \gamma = \const \in [0,1). \eqno(5.13)
$$

 \ It follows immediately from the direct definition of the norm for GLS

$$
\beta( x,y; \ p) \le \upsilon(p) \cdot r(x,y).
$$

 \ As long as

$$
N(X, \beta_p, \epsilon) \le N(X,\upsilon(p) r, \epsilon) = N(X, r, \epsilon/\upsilon(p))  \le C(N,r)  \ \epsilon^{-\gamma} \ \upsilon^{\gamma}(p),
$$
we deduce on the basis of theorem 3.1 for all the admissible values $  p  $

 $$
  Q(X,r,\xi;u) \le \frac{2 C(N,r)}{C(\lambda)} \ \upsilon^{\gamma}(p) \ W^{p + 1}(\gamma) \ u^{-p}, \ u \ge 1.
 $$

\vspace{4mm}

{\bf Theorem 5.2.}  We assert under the notations and conditions of this section

 $$
  Q(X,r,\xi;u) \le \frac{2 \ C(N,r)}{C(\lambda)}  \cdot \inf_{p \in [1,b)}  \left\{ \ \upsilon^{\gamma}(p) \ W^{p + 1}(\gamma) \ \ u^{-p} \right\},  \ u \ge 1.
\eqno(5.14)
 $$

\vspace{3mm}

 \ If in addition

$$
\lim_{h \to 0+} \sigma[r](h)  = 0,
$$
then

$$
{\bf P}(\kappa[\xi](h) > u) \le \ \frac{2 C(N,r)}{C(\lambda)} \  \times
$$

$$
\inf_{p \in [1,b)} \left[ \ \upsilon^{\gamma}(p) \ W^{p + 1}(\gamma) \ u^{-p}  \cdot \sigma[r](2 h)  \right], \ u \ge 1, \eqno(5.15)
$$
and as a consequence:  $  {\bf P} \left(\xi(\cdot) \in D[0,1]^d \right) = 1.  $\par

\vspace{3mm}

{\bf Example 5.2.}   Suppose the function  $ \ \upsilon(p), \ p \in [1,\infty) \ $ is such that

$$
\ \upsilon^{\gamma}(p) \le \exp \left(  C(1) \ p^{ 1 + v  } \right)
$$
for certain positive finite constant $ \ v. \ $ It is easy to calculate

$$
 Q(X,r,\xi;u) \le \exp \left( -C_2(v,\gamma,p) \ (\ln u)^{ 1 + 1/v }   \right), \ u \ge e. \eqno(5.16)
$$

\vspace{3mm}

{\bf Example 5.3.}   Suppose the distance function  $ \  r = r(x,y), \ x,y \in X \ $ is such that

$$
r = r(x,y) \asymp \sum_{j=1}^d |x(j) - y(j)|^{\alpha(j)}, \ \alpha(j) = \const \in (0,\infty).
$$

 \ The constant  $  \gamma $ may be calculated  by the formula

$$
\frac{1}{\gamma} =    \frac{1}{\alpha(1)}  +  \frac{1}{\alpha(2)}  + \ldots +   \frac{1}{\alpha(d)}. \eqno(5.17)
$$

 \ If for instance $ \ \forall j \ \alpha(j) = \alpha \in (0,\infty),   $ then $ \gamma = d/\alpha; \ $   and the used before conditions

$$
\lim_{h \to 0+} h^{-d} \sigma[r](h) = 0 \ \Leftrightarrow \gamma < 1
$$
are  equivalent to the following restriction:  $  \ \alpha > d. $ \par

 \vspace{3mm}

\section{ Auxiliary estimates III. Mixed moment estimates for normalized sums of independent random variables. }

 \vspace{3mm}

 \  We intend to extend in this section the results of penultimate one to the normalized sums of independent random variables. \par

 \vspace{3mm}

 \ We agree  to take in the sequel during this and next sections  that $  p(j) \ge 2; \  $ therefore $  b(j) \in (2, \infty]. $ Recall here the
famous {\it  Rosenthal's }  inequality, see [41],[22],[24],[33],[28] etc. \par

 \ Let $  \iota_i, \ i = 1,2,\ldots; \ \iota = \iota_1 $ be a sequence of the centered i., i.d. random variables with finite $ p^{th} $ moment.
 Rosenthal's inequality tell us:

$$
\sup_n \left| n^{-1/2} \ \sum_{i=1}^n \iota_i \ \right|_p \le C(R) \ \frac{p}{\ln p} \ |\iota|_p, \eqno(6.0)
$$
where $  C(R) $ is an absolute constant. Denote for brevity $   K(p) = K_R(p) := C(R) \ p/\ln p -   $ the so-called Rosenthal constant, more exactly,
Rosenthal's function on $ p. $\par

\ We  introduce in this connection  the following Rosenthal transformation  $  \psi_R(p)   $ for each $ \Psi - $ function  $  \ \psi(\cdot): $

$$
\psi_R(p) \stackrel{def}{=}   K_R(p) \cdot \psi(p) =  K_R[\psi](p):=  [C(R) \ p/\ln p ] \cdot \psi(p). \eqno(6.1)
$$

 \ It is clear that if $ \supp \psi(\cdot) < \infty, $ then  $ \ \psi_R(p) \asymp \psi(p). $\par

 \ The Rosenthal's inequality may be rewritten by means of our notations as follows

$$
\sup_n  || \ \left| n^{-1/2} \ \sum_{i=1}^n \iota_i \ \right| \ ||G\psi_R \le || \ \iota \ ||G\psi. \eqno(6.2)
$$
 \ Note that the last estimate in essentially non - improvable. \par

\vspace{3mm}

 \ We need more some generalization of Rosenthal's inequality on the multivariate mixed moments, alike in the fourth section.  Indeed,
let $ \ \{ \xi_i^{(j)}  \}, \ i = 1,2, \ldots; \ j = 1,2, \ldots, k \  $ be $  \ k \  $ tuple of infinite sequences of {\it centered} random variables.
We impose the following conditions of its independence: \\

\vspace{3mm}

 {\it for all the (fixed) upper index}  $ \  j \  $  {\it  the centered random variables } $  \ \{ \xi_i^{(j)}  \}, \ i = 1,2, \ldots; \hspace{4mm} \xi^{(j)} := \xi_1^{(j)} \ $
{\it are  independent and identical distributed. } \par

\vspace{3mm}

 \ Denote for every such  the index $  j  $

$$
S_n(j) = n^{-1/2} \sum_{i=1}^n \xi_i^{(j)}, \eqno(6.3)
$$
so that  by virtue of Rosenthal's inequality

$$
\sup_n  \left| \ S_n(j) \ \right|_p  \le K_R(p) \  \left| \ \xi^{(j)} \  \right|_p =  K(p) \  \left| \ \xi_1^{(j)}  \  \right|_p, \eqno(6.4)
$$
or equally

$$
\sup_n  || \  S_n(j) \ || G\psi_{j,R}  \le   || \ \xi^{(j)} \  || G\psi_j; \eqno(6.4a)
$$
and introduce  the following {\it mixed}  moment

$$
\nu = \nu \left(  \{  \xi_i^{(j)} \}, \{ p_j \}  \right) =  \nu \left( \ \vec{\xi}, \  \vec{p} \ \right) =
$$

$$
\nu(p_1,p_2,\ldots,p_k ) \stackrel{def}{=}  \sup_n {\bf E} \left[ \ \prod_{j=1}^k  \left|S_n (j) \right|^{p_j} \ \right], \ p_j = \const  \ge 2. \eqno(6.5)
$$

 \ Let $  \{ \ a(1), a(2), \ldots, a(k)  \ \}  $ be again the $ \ k - \ $ tuple of real numbers  greatest that one:  $  a(j) > 1, $
and such that  as before

$$
\sum_{j=1}^k \frac{1}{a(j)} = 1.
$$
 \ The set all of such the $ \ k \ $  tuples we denoted by $  A = A(k). $   We apply  again H\"older's inequality

$$
\nu \le \prod_{j=1}^k  \ \left\{ K_R[\psi]( a(j) \ p(j) )  \  \left| \xi_i^{j)} \ \right|_{a(j)p(j)} \right\}^{p(j)} =
$$

$$
 \prod_{j=1}^k  \left[\psi_{j,R} (a(j) \ p(j)) \right]^{p(j)} := Z(\vec{a}, \vec{p}),
$$
thus

$$
\mu \left( \ \vec{\xi}, \  \vec{p} \ \right) \le \inf_{\vec{a} \in A(k)} Z(\vec{a}, \vec{p}) =: \overline{Z}(\vec{p}). \eqno(6.6)
$$

 \ It follows immediately from the last inequality an uniform tail estimate

$$
 \sup_n {\bf P} \left(  \min_j  | S_n(j) | > u   \right) \le \frac{ \overline{Z}( \vec{p})}{u^{\sum_j p(j)}}, \ u > 0,  \eqno(6.7)
$$
if of course the right-hand side  of  (6.7) is not only is finite but also tends to zero as $   \ u \to \infty. $ \par

 \ As  a slight consequence:

$$
 \sup_n {\bf P} \left(  \min_j  | S_n(j) | > u   \right) \le  \ \inf_{\vec{p} \ge \vec{1}} \ \left[ \ \frac{ \overline{Z}( \vec{p})}{u^{\sum_j p(j)}} \ \right], \ u > 0.
  \eqno(6.8)
$$

\vspace{3mm}

\section{ Central Limit Theorem in the multivariate Prokhorov-Skorokhod space.}

\vspace{3mm}

 \ Let now $  \xi_i = \xi_i(x), \ \xi(x) = \xi_1(x), \ x \in X = [0,1]^d $ be a sequence of separable  stochastic continuous
independent  identically distributed (i; i.d)  {\it centered (mean zero): \ }  $ \ {\bf E} \xi_i(x) = 0, \ x \in X \ $ random fields with finite variance
at each the point: $ \ \forall x \in X \ \Var \xi_i(x) < \infty, \ \xi(x) := \xi_1(x). $  Denote

$$
S_n(x) = n^{-1/2} \sum_{i=1}^n \xi_i(x). \eqno(7.0)
$$

 \ Denote also by  $ \ S(x) = S_{\infty}(x)  \ $ the centered separable stochastic continuous Gaussian random field with at the same covariation
function as the source r.f.  $  \xi(x): $

$$
R(x,y) :=  {\bf E} S(x)  S(y) = {\bf E}  \xi(x)  \xi(y),  \ x,y \in X.
$$

 \ However, the stochastic continuity of the limiting Gaussian field $  \  S = S(x) \ $ is quite  equivalent to the continuity of its  covariation
function  $ \  R(x,y), \ $ which implies in turn the stochastic continuity of each r.f. $ \ \xi_i(x). $\par
 \  Note in addition

$$
{\bf E} S_n(x)  S_n(y) = R(x,y); \ \hspace{4mm} R(x,x) = \Var ( \xi_i(x)) = \Var ( S_n(x)), \ x \in X.
$$

 \ Evidently,  the finite-dimensional distributions of the sequence of r.f.  $  \  \{ \ S_n(\cdot) \} \ $ converge to ones for $ S(\cdot). $  Therefore, we
need to ground only the weak compactness of the sequence of distributions $  \  \Law \{ \ S_n(\cdot) \} \  $ in the Prokhorov-Skorokhod space,
formulated as before in the natural terms for the source r.f. $  \ \xi(x). \ $ More precisely:

\vspace{4mm}

{\bf Definition 7.1.} \\

\vspace{4mm}

 \ We will say as ordinary that the sequence $ S_n(\cdot), $ or more simple the alone r.f. $  \xi(x) $ satisfies CLT in the Prokhorov-Skorokhod space
$  D[0,1]^d  = D(X),  $  if  $ \xi(\cdot) \in D(X) $ with probability one and the sequence $ S_n(\cdot)  $ converges weakly (in distribution) to the
r.f. $  S(x). $ Briefly: $ \ \xi = \xi(x) \in CLT(PS) $ or more detail

$$
n \to \infty \ \Rightarrow S_n(\cdot) \stackrel{d}{\to} S_{\infty}(\cdot) = S(\cdot). \eqno(7.1)
$$

\vspace{4mm}

 \ The first examples of  $  \ CLT(PS) \ $  with statistical applications  belong to Yu.V.Prokhorov [40] and A.V.Skorokhod [43].
Another applications in the Monte-Carlo method and in the reliability theory may be found in [20]. The multivariate case $ d \ge 2 $ is
investigated in the famous article  of P.J.Bickel and M.J.Wichura  [3]. \par
 \ Denote

$$
\Delta_n[\xi](M)  =  \Delta_n[\xi](M) \left(x^{(2)} \right):= \left( S_n (x^{(2)}) - S_n(z_M)  \right), \ M \subset(T),
$$

$$
\Delta^{(i)}[\xi](M) = \Delta_n^{(i)}[\xi](M) \left(x^{(2)} \right):=   \left( \xi_i(x^{(2)}) - \xi_i(z_M)  \right);
$$

$$
 U[M](p) =  U[M; x^{(1)}, x^{(3)}](p)  := \sup_{x^{(2)} \in X} |\Delta[\xi](M)|_p, \ p \ge 2,
$$

$$
V[M](p) = V[M; x^{(1)}, x^{(3)}](p) := \sup_{x^{(2)} \in X} \sup_n |\Delta_n[\xi](M)|_p, \ p \ge 2,
$$

 \ We have

$$
\Delta_n[\xi](M)  :=   n^{-1/2} \sum_{i=1}^n \Delta^{(i)}[\xi](M),
$$
and we derive following  by virtue of Rosenthal's inequality

$$
V[M](p)  \le K_R(p) \cdot U[M](p) =: KU[M](p).  \eqno(7.2)
$$

 \ Let us introduce the following quasy-distance function, more exactly, the family of  quasy-distance functions

$$
 \delta \left( \ x^{(1)}, \ x^{(3)} \ \right) =  \delta \left(x^{(1)}, \ x^{(3)}; \vec{s} \right) =
\delta[\xi] \left(x^{(1)}, \ x^{(3)}; \vec{s} \right) \stackrel{def}{=}
$$

$$
\sup_{x^{(2)} \in X} {\bf E} \sup_n  \prod_{M \subset T} \left[ \  \left|  \ \Delta_n[\xi](M) \right|^{s(M)} \ \right].   \eqno(7.3)
$$

\vspace{4mm}

 \ This function allows a very simple estimate  by means of H\"older's inequality. Namely, if we denote

$$
\tilde{\delta} \left(x^{(1)}, \ x^{(3)}; \vec{s} \right) = \tilde{\delta}[\xi] \left(x^{(1)}, \ x^{(3)}; \vec{s} \right)
 \stackrel{def}{=}
$$

$$
 \sup_{x^{(2)} \in X} \prod_{M \subset T} \left[ KU[M]  (\alpha(M) \ s(M) )  \right]^{1/\alpha(M)},   \eqno(7.4)
$$
then

$$
 \delta[\xi] \left(x^{(1)}, \ x^{(3)}; \vec{s} \right)  \le  \tilde{\delta} \left(x^{(1)}, \ x^{(3)}; \vec{s} \right).   \eqno(7.5)
$$

 \ Further,

$$
 \delta \left(x^{(1)}, \ x^{(3)}; \vec{s} \right)  \le   \tilde{\delta} \left(x^{(1)}, \ x^{(3)}; \vec{s} \right) \le
 \delta^+  \left(x^{(1)}, \ x^{(3)}; \vec{s} \right),
$$
where

$$
\delta^+  \left(x^{(1)}, \ x^{(3)}; \vec{s} \right) = \delta^+[\xi]  \left(x^{(1)}, \ x^{(3)}; \vec{s} \right) \stackrel{def}{=}
$$

$$
\inf_{\vec{\alpha} \in A(T)} \sup_{x^{(2)} \in X} \prod_{M \subset T}  \left[ KU[M]  (\alpha(M) \ s(M) )  \right]^{1/\alpha(M)},  \eqno(7.6)
$$

where  as was descripted $  \ M \subset T, $
$$
\vec{\alpha} \in A(T) \ \Leftrightarrow \   \alpha(M) > 1, \ \sum_M 1/\alpha(M) = 1.
$$

\vspace{4mm}

 \ Assumptions and notations:  $ \  p := \sum_M s(M) > 0, \  $

$$
N(X,\delta[\xi](\cdot, \cdot; \vec{s}), \ \epsilon) \le C_{\delta}(N) \ \epsilon^{-\gamma}, \ \epsilon \in (0,1), \
$$

$$
 \gamma = \const  \in (0,1), \ C_{\delta}(N)  = \const \in (0,\infty),   \eqno(7.7)
$$
following,
$$
\lambda(u) = C_{\delta}(\lambda) \  u^{\sum_M s(M) }= C_{\delta}(\lambda) \  u^p ,  \ u > 0;  \ C(\lambda) =  \const \in (0,\infty).  \eqno(7.8)
$$

\vspace{4mm}

{\bf Theorem 7.1.} We get under our notations and conditions (7.7) etc.  by virtue of theorem 2.2

 $$
\sup_n  Q(X,\delta, S_n; \ u) \le \frac{2 C_{\delta}(N)}{C_{\delta}(\lambda)} \ W^{p + 1}(\gamma) \ u^{-p}, \ u \ge 1.   \eqno(7.9)
 $$

 \ If in addition

$$
\lim_{h \to 0+} \sigma[\delta](h)  = 0,   \eqno(7.10)
$$
then

$$
\sup_n {\bf P}(\kappa[S_n](h) > u) \le \
$$

$$
\frac{2 C(N)}{C(\lambda)} \ W^{p + 1}(\gamma) \ u^{-p}  \cdot \sigma[\beta](2 h), \ u \ge 1,  \eqno(7.11)
$$
and as a consequence the r.f. $  \ \xi(\cdot) $ satisfies the CLT in the Prokhorov-Skorokhod  space $  \ X = D[0,1]^d. $ \par

\vspace{4mm}

{\bf Remark 7.1}.  Obviously, instead the condition (7.7) may be used the following more easily verified one

$$
N(X, \tilde{\delta}[\xi](\cdot, \cdot; \vec{s}), \ \epsilon) \le C_{\tilde{\delta}}(N) \ \epsilon^{-\gamma}, \
$$

$$
\epsilon \in (0,1), \ \gamma = \const  \in (0,1), \ C_{\tilde{\delta}}(N)  = \const \in (0,\infty),   \eqno(7.12)
$$
or in turn

$$
N(X, \delta^+[\xi](\cdot, \cdot; \vec{s}), \ \epsilon) \le C_{\delta^+}(N) \ \epsilon^{-\gamma}, \
$$

$$
 \epsilon \in (0,1), \ \gamma = \const  \in (0,1), \ C_{\delta^+}(N)  = \const \in (0,\infty).   \eqno(7.13)
$$

\vspace{4mm}

{\bf Remark 7.2.} The conditions of theorem 7.1 are satisfied for instance for the so-called empirical random fields, see [40], [43];
 see also  [42];  as well as for multiple parametric integral computation by means of Monte-Carlo method  [20]. \par

 \vspace{3mm}

\section{ Concluding remarks.  }

 \vspace{3mm}

{\bf A. \ Necessary and sufficient condition for the  Prokhorov-Skorokhod continuity of random fields.} \par

 \vspace{3mm}

 \ Let $  \ \xi = \xi(x), \ x \in X $ be stochastic continuous  separable r.f.  The classical module of (uniform) continuity
$ \ \omega[\xi](h), \ h \in [0,1] \ $ for the r.f. $  \ \xi(\cdot) \ $  is defined as follows.

$$
\omega[\xi](h) \stackrel{def}{=} \sup_{x,y: |x-y| \le h} |\xi(x) - \xi(y)|. \eqno(8.0)
$$

 \ The r.f. $ \ \xi(\cdot) $ is continuous with probability one, or equally has a continuous almost surely sample path, iff

$$
{\bf P} \left(\lim_{h \to 0+} \omega[\xi](h) = 0  \right) = 1.  \eqno(8.1)
$$
 \ The last relation is  completely equivalent to the following  natural condition

$$
\lim_{h \to 0+} {\bf E} \arctan \omega[\xi](h) = 0, \eqno(8.2)
$$
as long as  $ \ \omega[\xi](h), \ h \in [0,1] $ is monotonically increasing function. \par
 \ Of course,   our set $  \ [0,1]^d \ $ may be replaced by arbitrary complete compact metric space. \par

 \ Analogous statement, with at the same proof, holds true for the Prokhorov-Skorokhod space.\par

\vspace{4mm}

{\bf Proposition 8.1.} Let  $ xi(x), \ x \in X   $ be arbitrary separable continuous in probability r.f. In order to it belongs to the
Prokhorov-Skorokhod space $ \ D[0,1]^d \ $ almost everywhere, is necessary and sufficient

$$
\lim_{h \to 0+} {\bf E} \arctan \kappa[\xi](h) = 0.  \eqno(8.3)
$$

 \ Analogously may be formulated the criterion for weak compactness the family of distributions of stochastic continuous r.f. $  \{  \xi_n(x) \}, \ x \in X $
in the  Prokhorov-Skorokhod space $ \ D[0,1]^d \ $ with converges all the finite-dimensional distributions:

$$
\lim_{h \to 0+}  \sup_n {\bf E} \arctan \kappa[\xi_n](h) = 0.  \eqno(8.3)
$$

\vspace{4mm}

{\bf B. \ Factorability of module of continuity.   } \par

 \ It is known,[28], chapter 4, that for arbitrary continuous with probability one random field $  \ \xi(x), \ x \in X  $ its
classical (uniform) module of continuity  $  \ \omega[\xi](h), \ h \in [0,1] $  allows a {\it factorization}:  there exists  a random variable
$ \ \theta \  $ and  {\it a  non-random } module  of continuity $ \ \omega_o = \omega_o(h), \ 0 \le h \le 1,  $  such that

$$
\omega[\xi](h) \le \theta \cdot  \omega_o(h).   \eqno(8.4)
$$
 \ Of course, the converse conclusion is also true. \par

 \ Analogous factorization, with at the same proof,  remains true for the Prokhorov-Skorokhod space. In detail,  the  r.f.  $ \ \xi(x) \ $
belongs to the space $  \ D[0,1]^d $ a.e., or equally

$$
\lim_{h \to 0+}  \kappa[\xi])(h) = 0  \ (\mod P),
$$
if and only if

$$
\kappa[\xi](h) \le \Theta \cdot \omega_0(h)
$$
for some finite r.v. $  \ \Theta \  $ and for certain non-random module of continuity $  \ \omega_0(h). $\par

 \vspace{3mm}
 {\bf C.  Minimum estimates by means of $  B(\phi)  $ spaces. } \par
\vspace{3mm}

 \ The  {\it exponential} tail estimates for distribution for minimum of the finite set of random variables
based on the theory of so-called   $  B(\phi)  $ spaces with an application to the theory of discontinuous random fields
may be found in [28], chapter 3, section 3.16; [29],[36].\par

 \vspace{3mm}
 {\bf D.  Compact embedded support.} \par
\vspace{3mm}

 \ Let $  \ B \ $ be separable Banach space and let $  \  \mu \ $ be sigma-finite Borelian measure on the space $  B. $ It is known, see [28], chapter 4; [31]
  that there exists a separable compact embedded Banach subspace $ \ Q \ $ of the space $  B  $ which  is complete support of the measure $ \ \mu:   $\par

$$
\mu(B \setminus Q) = 0.
$$
  \ Recall that the Banach subspace $  \ Q \ $ is said to be  compact embedded into Banach space $ \ B, \ $ if any bounded set of the space $  \ Q  \ $
is pre-compact set of the space $  \ B. \ $ \par
 \ As for the Linear Topological Spaces: note that this  proposition is not true still for the classical Schwartz space $ \ S \ $  of  infinite differentiable
functions having a finite support. \par
 \ It is more interesting to note that this statement  holds true yet for the considered in this report Prokhorov-Skorokhod spaces $  \ D[0,1]^d, \ $
despite the ones are not  Linear Topological Spaces! \par

 \vspace{4mm}

 \ \ \ \ \ \ \ \ \ \ \ \ \ \ \ \ \ \ \ \ \ \ \ \ \ \ \ \ \ \ {\bf References.}

 \vspace{4mm}

{\sc 1. Bennet C., Sharpley R.}  {\it  Interpolation of operators.} Orlando, Academic
Press Inc., (1988).\\

{\sc  2. P.H.Bezandry and X.Fernique.} {\it Sur la propriete de la limite centrale dans D[0; 1].}
 Ann. Inst. Henri Poincare, 28, (1992), no. 1, 31-46.\\

{\sc 3. P.J.Bickel and M.J. Wichura.} (1971). {\it Convergence criteria for multiparameter stochastic processes and some applications.}
 Ann. Math. Statist., 42, 1656-1670. \\

{\sc 4. P. Billingsley.} (1968). {\it Convergence of Probability Measures.} New York, John Wiley and Sons.\\

{\sc 5. P. Billingsley.} (1971). {\it Weak Convergence of Measures: Applications in Probability.}  Philadelphia. SIAM, New York-London. \\

{\sc 6. M.Bloznelis and V.Paulauskas.}  {\it On the central limit theorem in D[0; 1].}
Statistics and Probability Letters, 17, (1993), 105 -111.\\

{\sc 7. M.Bloznelis and V.Paulauskas.}  {\it  Central limit theorem in Skorokhod spaces
and asymptotic strength distribution of fiber bundles.} (1994), Sixth Intern. Vilnius Conference on Probab. Th. Math.Statist. Eds. B.Grigelionis, J.Kubilius,
H.Pragarauskas, V.Statulevichius, VSP/TEV, Utrecht/Vilnius, 75-87.\\

{\sc 8. M.Bloznelis and V.Paulauskas. }    {\it Central limit theorem in D[0,1].} Internet
Publication, 2015, PDF.\\

{\sc 9. Buldygin V.V., Kozachenko Yu.V. }  {\it Metric Characterization of Random
Variables and Random Processes.} 1998, Translations of Mathematics Monograph, AMS, v.188.\\

{\sc 10. N.N. Chentsov.}  {\it Weak convergence of stochastic processes whose trajectories
have no discontinuities of the second order.} Theory Probability Applications., 1-3, (1956),  140-143.\\

{\sc 11. H.E.Daniels.}  {\it The statistical theory of the strength of bundles of threads.} Proc.
Roy. Society. London, A, 183, (1945), 404-435.\\

{\sc 12. Dudley R.M.}  {\it Uniform Central Limit Theorem.} Cambridge University Press, (1999).\\

{\sc 13. Fernique X. (1975).}  {\it Regularite des trajectoires des function aleatiores gaussiennes. } Ecole de Probablite de Saint-Flour, IV-1974,
Lecture Notes in Mathematic. 480, 1-96, Springer Verlag, Berlin.\\

{\sc 14. Fernique X.}  {\it Regularite de fonctions aleatoires non gaussiennes.} Ecolee de Ete
de Probabilite de Saint-Flour XI-1981. Lecture Notes in Math., 976, (1983), 174, Springer, Berlin.\\

{\sc 15. Fernique X.} {\it Caracterisation de processus de trajectoires majores ou continues.}
 Seminaire de Probabilit?s XII. Lecture Notes in Math. 649, (1978), 691-706,  Springer, Berlin.\\

{\sc 16. A. Fiorenza.}   {\it Duality and reflexivity in grand Lebesgue spaces. } Collect. Math.
{\bf 51,}  (2000), 131-148.\\

{\sc 17. A. Fiorenza and G.E. Karadzhov.} {\it Grand and small Lebesgue spaces and
their analogs.} Consiglio Nationale Delle Ricerche, Instituto per le Applicazioni
del Calcoto Mauro Picone”, Sezione di Napoli, Rapporto tecnico 272/03, (2005).\\

{\sc 18. Frolov A.S., Chentzov N.N.} {\it On the calculation by the Monte-Carlo
method definite integrals depending on the parameters.} Journal of Computational Mathematics and Mathematical Physics, (1962), V. 2,
Issue 4, p. 714-718 (in Russian).\\

{\sc 19. I.I. Gikhman and A.V. Skorokhod. }  (1965).  {\it Introduction to the theory of random processes, } (Russian). Moscow, Nauka, 654 pp. English edition: \\
{\sc I.I.Gikhman and A.V. Skorokhod,} (1969), {\it Introduction to the theory of random
processes,} Philadelphia, W.B. Saunders Co., xiii+516 pp.\\

{\sc 20. Grigorjeva M.L., Ostrovsky E.I.}  {\it Calculation of Integrals on discontinuous functions by means of depending trials method.}
Journal of Computational Mathematics and Mathematical Physics, (1996), V. 36, Issue 12, p. 28-39 (in Russian).\\

{\sc 21. Daniel J. Greenhoe.} {\it Properties of distance spaces with power triangle inequalities.}
 arXiv:1610.07594 math.CA; 26.10.2016.\\

{\sc 22.  R.Ibragimov and Sh. Sharakhmetov.}  {\it On an exact constant for the Rosenthal
inequality. } Theory Probability Applications., 42, (1997), pp. 294-302.\\

{\sc 23. T.Iwaniec and C. Sbordone.} {\it On the integrability of the Jacobian under minimal
hypotheses. } Arch. Rat.Mech. Anal., 119, (1992), 129-143.\\

{\sc 24. W. B. Johnson, G. Schechtman, and J. Zinn. } {\it Best constants in moment
inequalities for linear combinations of independent and exchangeable random
 variables. } Ann. Probability, 13, (1985), pp. 234- 253.\\

{\sc 25. Kolmogorov A.N. }  {\it On the Skorokhod convergence.} Theory of probability and
its applications. V.1, (1956), pp. 239-247 (Russian), pp. 215-222, (English).\\

{\sc 26. Kozachenko Yu. V., Ostrovsky E.I. }  (1985). {\it The Banach Spaces of random Variables of subgaussian Type. } Theory of Probability.
and Math. Stat. (in Russian). Kiev, KSU, 32, 43-57.\\

{\sc 27. Petr Lachout.} {\it Billingsley-type tightness criteria for multi - parameter stochastic processes.}
Cybernetics, Vol. 24 (1988), No. 5, 363-371. \\

{\sc 28. Ostrovsky E.I. } (1999). {\it Exponential estimations for Random Fields and its
applications,} (in Russian). Moscow-Obninsk, OINPE.\\

{\sc 29. Ostrovsky E. and Sirota L.} {\it Vector rearrangement invariant Banach spaces
of random variables with exponential decreasing tails of distributions.} \\
 arXiv:1510.04182v1 [math.PR] 14 Oct 2015 \\

{\sc 30. Ostrovsky E. and Sirota L.}  {\it Non-asymptotical sharp exponential estimates
for maximum distribution of discontinuous random fields. } \\
 arXiv:1510.08945v1 [math.PR] 30 Oct 2015 \\

{\sc 31. Ostrovsky E.I.}  {\it About supports of probability measures in separable Banach
spaces.} Soviet Math., Doklady, (1980), V. 255, $ \ N^0 \ $ 6, p. 836-838, (in Russian).\\

{\sc 32. Ostrovsky E. and Sirota L.} {\it Criterion for convergence almost everywhere,
with applications.} \\
arXiv:1507.04020v1 [math.FA] 14 Jul 2015.\\

{\sc 33. Ostrovsky E. and Sirota L.}  {\it Schl\"omilch and Bell series for Bessel's functions, with probabilistic applications.} \\
 arXiv:0804.0089v1 [math.CV] 1 Apr 2008\\

{\sc 34. Ostrovsky E. and Sirota L. } {\it Sharp moment estimates for polynomial martingales. } \\
arXiv:1410.0739v1 [math.PR] 3 Oct 2014\\

{\sc 35. Ostrovsky E., Rogover E. } {\it Exact exponential bounds for the random field
maximum distribution via the majorizing measures (generic chaining).} \\
 arXiv:0802.0349v1 [math.PR] 4 Feb 2008 \\

{\sc 36. Ostrovsky E. and Sirota L. } {\it   Entropy and Grand Lebesgue Spaces approach for the problem  of Prokhorov-Skorokhod continuity of
discontinuous random fields. }\\
arXiv:1512.01909v1 [math.Pr] 7 Dec 2015 \\

{\sc 37. M.Bloznelis,  V.Paulauskas.}  {\it On the central Limit Theorem for multi-parameter stochastic processes. } In: Probability in Banach spaces, {\bf 9,}
 Newhaus, Straf, (1971), 155-172.\\

{\sc 38. V.Paulauskas and Ch.Stieve.}  {\it On the central limit theorem in D[0; 1] and in
D([0; 1];H). } Lietuvos matem. rink., 30, (1990), 267-279.\\

{\sc 39. Pizier G. }  {\it Condition d' entropic assupant la continuite de certains processus et
applications a lanalyse harmonique.} Seminaire d analyse fonctionalle. (1980), Exp.13, p. 23-34. \\

{\sc 40. Prokhorov Yu. V.}  {\it Convergence of random processes and limit theorems in
probability.}  Theory of probability and its applications. V.1, (1956), pp. 177-238,
(Russian), pp. 151-214, (English).\\

{\sc  41. H. P. Rosenthal.} {\it On the subspaces of L(p) (p > 2) spanned by sequences of
independent random variables. } Israel J. Math., 8  (1970), pp. 273-303.\\

{\sc  42. Serik Sagitov.} {\it Weak Convergence of Probability Measures.}  Internet Publication, April 2015, PDF.
Chalmers University of Technology and Gothenburg University.\\

{\sc 43. Skorokhod A.V. } {\it Limit theorems for stochastic processes.} Theory Probability Applications,
 1956, t. 1, no. 3, p. 289-319 (in Russian); English translation.: Theory
Probability Applications., 1956, v. 1, no. 3, p. 261-290.\\

{\sc 44. Talagrand M.}  (1996). {\it Majorizing measure: The generic chaining.} Ann.
Probability., 24 1049-1103. MR1825156 \\

{\sc 45. Talagrand M.} (2001). Majorizing Measures without Measures. Ann. Probability.,
29, 411-417. MR1825156 \\

\end{document}